# OPTIMIZATION OF SYSTEM OF NONLINEAR SECOND ORDER DIFFERENTIAL INEQUALITIES

**Elimhan N. Mahmudov[1,2], S. Demir Saglam[3],**

[1]*Department of Mathematics and Engineering, Istanbul Technical University, Istanbul, Turkey,*

[2]*Institute of Control Systems, Azerbaijan National Academy of Sciences, Azerbaijan,*

[3]*Department of Mathematics, Istanbul University, Istanbul, Turkey.*

**Abstract**. This paper deals with the optimization of Bolza problem with a system of convex and nonconvex, discrete and differential state variable inequality constraints of second order by deriving necessary and sufficient conditions for optimality. According to proposed discretization method and equivalence theorems for subdifferential inclusions, the problem with a system of discrete-approximation inequalities are investigated which highly contributes to the derivation of adjoint discrete inclusions generated by given system of nonlinear inequality constraints. Moreover, in the limit case, we obtain sufficient conditions for optimality of the continuous problem. A numerical example is presented to illustrate the theoretical result.

**Key words:** Differential inequality, Inclusion, Adjoint, Approximation.

**1. Introduction and Needed Facts**

In this paper, a particular optimization problem is investigated whose constraints are defined by a system of second order discrete and differential nonlinear inequality constraints. In general, optimization problems with a system of inequality constraints have wide applications such as in economy, engineering, industry and game theory. Especially, differential inequalities appear naturally in several areas of applied mathematics, like in mathematical economics in problems of resource allocation and in the study of planning procedures, in mechanics in the study of elastoplastic systems and in differential games. Many mathematical problems in literature can be formulated as variational inequality problems with continuous multivalued mappings and their relationships with other general problems of nonlinear analysis [1,2,6,7,9,13-18,31,35,37-39].

One general approach to such optimization relies on the representation of the sublevel sets of the inequality constraints via set-valued mappings. This analysis reveals the hidden relationship between optimal control methods with inequality constraints and set-valued mappings; in the paper [5], sensitivity relations are obtained for the Mayer problem associated with the differential inclusion $x' \in F(x)$ and applied to derive optimality conditions. The authors first application concerns the maximum principle and consists in showing that a dual arc can be constructed for every element of the superdifferential of the final cost.

Actually, the tough part of the problems with higher order differential inequalities and/or inclusions is establishing higher order adjoint inclusions of Euler-Lagrange type and convenient transversality conditions [26-28]. For this reason, qualitative properties of second order differential inequalities and/or inclusions are predominantly considered in literature (see [3-5,10,12,21-23,25,33] and references therein). In the paper [12] nonlinear differential inequality is formulated, for whose solutions, the growth and decay rate is estimated. For the dynamical systems and nonlinear evolution equations in Banach spaces, the latter inequality is essential and applicable to investigation of global existence of solutions to nonlinear partial differential equations. In the paper [34], it is shown that if convexity is taken in the generalized sense, a differential characterization in the classical sense can be obtained. As a particular case of a general theorem concerning second order differential inequalities, a recent result of Tchaplygin, concerning linear differential inequalities is obtained. In the work [36], using a classical result on linear differential inequalities, some results for Riccati inequality and the second order linear differential inequality are established.

In the work [11], a multiobjective optimization problem with differentiable equality constraint functions and locally Lipschitz objective and inequality constraint functions. It is well known that modified Lagrangian functionals can be utilized in the method of Lagrange multipliers to solve the finite dimensional convex optimization problem. Lately in mechanics, infinite-dimensional variational inequalities are solved by this method. Some stability results for generalized vector quasi-variational inequality problems are derived in the paper [19] and the upper semi-continuity properties of the solution set are established for perturbed generalized vector quasi-variational inequality problems.

The paper [32] gives a second order overview of a significant class of variational systems in finite-dimensional and infinite-dimensional spaces, which is particularly important for the analysis of optimization and equilibrium problems with equilibrium constraints. In addition, the paper describes systems of this type through variational inequalities over polyhedral convex sets.

These optimization problems are extensions of works done in Mahmudov's earlier papers [20-29], where the main concerns are the necessary and sufficient conditions for optimality of the Bolza problem with higher order discrete and differential inclusions.

The stated problems and obtained conditions for optimality are novel. This paper consists of three parts.

In the first part, an optimal control problem is formulated where the system dynamics are given by a system of second order discrete inequalities (SSDSI).

In the second part of the paper, an optimization problem with a system of second order differential inequalities (SSDFI) involving first and second order derivatives of searched functions is formulated.

In the third part, a discrete-approximation problem for a SSDSI and sufficient conditions of optimality for a SSDFI are established.

Thus, for the reader's convenience, notions like local tent, convex upper approximation (CUA) of non-smooth functions, etc. from Mahmudov's book [21] and papers [22,24,26] are given in Section 1.

In Section 2, the Cauchy problems for a SSDSI and SSDFI are devised.

In Section 3, the optimality problem for the established SSDSI is reduced to a problem with finite number of geometric constraints. By means of convex and nonsmooth analysis, necessary and sufficient conditions for optimality of the SSDFI are found. Under non-degeneracy assumptions imposed on the optimal trajectory for a convex problem, necessary and sufficient conditions are proven. Consequently, under the existence of convex upper approximation and local tents, the optimality conditions for SSDSI are obtained.

In Section 4, difference operators of first and second order are utilized in the problem for a SSDFI and associated with the SSDSI. Apparently, particular equivalence theorems are needed to pass to the discrete-approximation problem, which in turn links the main results of systems of second order discrete and discrete-approximation problems. By using local tents and CUAs, similar results are proven for nonconvex problems. Obviously, this can also have a significant place in numerical methods.

In Section 5, we derive sufficient conditions for optimality of a SSDSI which are built upon the limit case for the optimality conditions of a SSDFI. Applying the second order adjoint inclusions generated by inequality constraints and the discrete approximations method, sufficient optimality conditions for SSDFI are obtained. Usually, the system of second order adjoint differential inclusion involves some auxiliary adjoint variables. It appears that in some concrete problems, the adjoint differential inclusion only involves the "main" variable, i.e. $x^*(\cdot)$.

In what follows, the necessary standard notions can be found in [20-29]. As usual, $\mathbb{R}^n$ is an $n$-dimensional Euclidean space, $\langle x, v \rangle$ is the inner product of elements $x, v \in \mathbb{R}^n$ and $(x, v)$ is a pair of $x, v$. By int $A$ and ri $A$, we denote the interior and relative interior of a set $A \subset \mathbb{R}^{3n}$.

**Definition 1.1**[21] A convex cone $K_A(z_0), z = (x, v_1, v_2)$ is called the cone of tangent directions at a point $z_0 = (x^0, v_1^0, v_2^0) \in A$ if from $\bar{z} = (\bar{x}, \bar{v}_1, \bar{v}_2) \in K_A(z_0)$ it follows that $\bar{z}$ is a tangent vector to the set $A$, i.e., there exists a function $\varphi(\lambda) \in \mathbb{R}^{3n}$ satisfying $z_0 + \lambda \bar{z} + \varphi(\lambda) \in A$ for sufficiently small $\lambda > 0$, where $\lambda^{-1} \varphi(\lambda) \to 0$, as $\lambda \downarrow 0$.

Evidently, for a convex set $A$ at a point $(x^0, v_1^0, v_2^0) \in A$ setting $\varphi(\lambda) \equiv 0$ we have
$$K_A(z_0) = \{(\bar{x}, \bar{u}, \bar{v}): \bar{x} = \lambda(x_1 - x^0), \bar{u} = \lambda(u_1 - v_1^0), \bar{v} = \lambda(v_1 - v_2^0), \lambda > 0\}, \forall (x_1, u_1, v_1) \in A.$$
As usually, $K_A^*(x^0, v_1^0, v_2^0)$ is the dual cone to a cone of tangent vectors $K_A(x^0, v_1^0, v_2^0)$.

**Definition 1.2** [21] For an arbitrary nonempty subset $P \subset \mathbb{R}^{3n}$, the cone defined by
$$\text{cone } P = \{(\bar{x}, \bar{v}_1, \bar{v}_2): \bar{x} = \alpha x, \bar{v}_1 = \alpha v_1, \bar{v}_2 = \alpha v_2, (x, v_1, v_2) \in P, \alpha > 0\}$$
is called the cone generated by $P$.

By Definition 1.1, we have already seen that the cone of tangent directions involve directions for each of which, there exists a function $\varphi(\lambda)$. However, this is not enough to determine the properties of the set $A$. Nonetheless, the next notion of a local tent will predetermine a mapping in $A$ for nearest tangent directions.

**Definition 1.3** [21] The cone $K_A(z_0)$ is called the local tent if for any $\bar{z}_0 \in \text{ri } K_A(z_0)$ there exists a convex cone $K \subseteq K_A(z_0)$ and a continuous mapping $\gamma(\bar{z})$ defined in the neighbourhood of the origin such that

(1) $\bar{z}_0 \in \text{ri } K$, $\text{Lin } K = \text{Lin } K_A(z_0)$, where $\text{Lin } K$ is the linear span of $K$,

(2) $\gamma(\bar{z}) = \bar{z} + r(\bar{z})$, $r(\bar{z}) \|\bar{z}\|^{-1} \to 0$ as $\bar{z} \to 0$,

(3) $z_0 + \gamma(\bar{z}) \in A$, $\bar{z} \in K \cap S_\varepsilon(0)$ for some $\varepsilon > 0$, where $S_\varepsilon(0)$ is the ball of radius $\varepsilon$.

**Definition 1.4** According to [21], $h(\bar{z}, z)$ is called a convex upper approximation (CUA) of the function $W: \mathbb{R}^{3n} \to \mathbb{R}^1 \cup \{\pm\infty\}$ at a point $z \in \text{dom} W = \{z: |W(z)| < +\infty\}$ if

(1) $h(\bar{z}, z) \geq V(\bar{z}, z)$ for all $\bar{z} \neq 0$,

(2) $h(\bar{z}, z)$ is a convex closed (or lower semicontinuous) positive homogeneous function in $\bar{z}$, where

$$V(\bar{z}, z) = \sup_{r(\cdot)} \limsup_{\lambda \downarrow 0} \frac{1}{\lambda} \left[ W(z + \lambda \bar{z} + r(\lambda)) - W(z) \right], \ \lambda^{-1} r(\lambda) \to 0.$$

Here, the exterior supremum is taken on all $r(\lambda)$ such that $\lambda^{-1} r(\lambda) \to 0$ as $\lambda \downarrow 0$.

**Definition 1.5** [21] The set defined as
$$\partial h(0, z) = \left\{ z^* \in \mathbb{R}^{3n} : h(\bar{z}, z) \geq <\bar{z}, z^*>, \bar{z} \in \mathbb{R}^{3n} \right\}$$
is called a subdifferential of the function $W$ at a point $z$ and is denoted by $\partial W(z)$. It should be noted that $h(\bar{z}, z)$ is the support function of $\partial W(z)$. Thus, $h(\bar{z}, z)$ and $\partial W(z)$ define each other one to one and the function $h(\cdot, z)$ defined by equation
$$h(\bar{z}, z) = \sup_{z^*} \left\{ \langle \bar{z}, z \rangle : z^* \in \partial h(0, z) \right\} = \sup_{z^*} \left\{ \langle \bar{z}, z \rangle : z^* \in \partial W(z) \right\}$$
must be a CUA of $W$ at $z$. Besides, if $h_1$ and $h_2$ are CUAs for the function $W$ at a point $z$ and $h_1 \geq h_2$ then $\partial_1 W(z) \supseteq \partial_2 W(z)$, where $\partial_1 W(z)$ and $\partial_2 W(z)$ are the subdifferentials defined by $h_1$ and $h_2$, respectively [20,21]. It is evident that $h_1$ is an approximation worse than $h_2$ in some neighborhood of $z$.

A function $W$ is called a proper function if it does not assume the value $-\infty$ and is not identically equal to $+\infty$. Due to its simpleness, this definition is very adventageous. It is to be noted [8,21,28] that this definition corresponds to the classical definition of a subdifferential in the convexity of $W$. Furthermore, a subdifferential can be constructed in various ways for different function classes. The primal notion of a subdifferential for general nonconvex functions was put forth by Clarke (see [8] and references therein) in nonsmooth analysis. He peculiarly established a broad subdifferential calculus for generalized gradients of locally Lipschitzian functions defined on Banach spaces. On the other hand, Mordukhovich's subdifferential [30] includes generalized differentials significant in pure and applied analysis.

## 2. Problem Statement

In this section, we formulate the so-called problems SSDSI and SSDFI. In the beginning, we consider the following second order discrete model denoted by (PD):

$$\text{minimize } \sum_{t=2}^{N-1} f(x_t, t), \tag{1}$$

(PD) 
$$W_k(x_t, x_{t+1}, x_{t+2}) \leq 0, \ t = 0, \ldots, N-2, \ k \in I = \{1, \ldots, m\}, \tag{2}$$

$$x_0 = \bar{\upsilon}_0, \ x_1 = \bar{\upsilon}_1, \tag{3}$$

where $x_t \in \mathbb{R}^n$, $f(\cdot, t)$ ($t = 0, \ldots, N-2$) are real-valued continuously differentiable convex functions, $f(\cdot, t): \mathbb{R}^n \to \mathbb{R}^1 \cup \{\pm\infty\}$, $W_k : \mathbb{R}^n \to \mathbb{R}^1 (k \in I)$ and $f(\cdot, N-1)$ are continuous convex functions, $N$ is a fixed natural number, and $\bar{\upsilon}_0, \bar{\upsilon}_1$ are fixed vectors.

In fact, the condition (3) is the discrete analogy of Cauchy initial conditions for a SSDSI. A sequence $\{x_t\}_{t=0}^N = \{x_t : t = 0, 1, \ldots, N\}$ is called a feasible trajectory of the given problem (1) – (3). The question is to find a feasible trajectory of problem (1)-(3) which minimizes $\sum_{t=2}^{N-1} f(x_t, t)$.

The problem (1) – (3) is a convex if $W_k (k \in I)$ and $f(\cdot, N-1)$ are continuous convex functions and $f(\cdot, t)(t = 2, \ldots, N-2)$ is continuously differentiable proper convex function.

**Condition** $H_1$ Suppose that in the convex problem (PD) for point $x_t^0 \in \mathbb{R}^n$, one of the following is satisfied:
(1) $(x_t^0, x_{t+1}^0, x_{t+2}^0) \in \text{ri } F$, $x_t^0 \in \text{ri }(\text{dom } f(\cdot,t))$, $t = 0,...,N-2$,
(2) $(x_t^0, x_{t+1}^0, x_{t+2}^0) \in \text{int } F$, $F = \{(x, v_1, v_2): W_k(x, v_1, v_2) \leq 0, k \in I\}$ (with the possible exception of one fixed $t_0$) and $f(\cdot, t)$ are continuous at $x_t^0$.

In what follows, this condition will be called the non-degeneracy condition.

**Condition** $H_2$ Assume that in the nonconvex case of the problem (PD) the cones of tangent directions $K_F(\tilde{x}_t, \tilde{x}_{t+1}, \tilde{x}_{t+2})$ are local tents, where $\tilde{x}_t$ are the points of the optimal trajectory $\{\tilde{x}_t\}_{t=0}^N$. In addition, suppose that the functions $W_k(x_t, x_{t+1}, x_{t+2}), k \in I$ and $f(\cdot, N-1)$ admit a continuous convex upper approximations (CUAs) $h_{W_k}^t(\cdot, \tilde{x}_t, \tilde{x}_{t+1}, \tilde{x}_{t+2})$ and $h_{N-1}(\cdot, \tilde{x}_{N-1})$, respectively [21] which ensure that the subdifferentials $\partial_{(x,v_1,v_2)} W_k(\tilde{x}_t, \tilde{x}_{t+1}, \tilde{x}_{t+2})$ $= \partial h_{W_k}^t(0, \tilde{x}_t, \tilde{x}_{t+1}, \tilde{x}_{t+2})$ and $\partial f(\tilde{x}_{N-1}, N-1) = \partial h_{N-1}(0, \tilde{x}_{N-1})$ are defined. Here, it is accepted that $W_k(\tilde{x}_t, \tilde{x}_{t+1}, \tilde{x}_{t+2}) = 0$. Further, let $f(\cdot, t)$ ($t = 2,...,N-2$) be continuously differentiable nonconvex function, that is $h_t(\overline{x}, \tilde{x}_t) = \langle \overline{x}, f'(\tilde{x}_t, t) \rangle$.

Besides, if $W_k, k \in I$ is a continuously differentiable nonconvex function, then the cone $K_F(z_0)$, $z_0 \in F$ defined as
$$K_F(z_0) = \{\overline{z}: \langle \overline{z}, W'_{kz}(z_0) \rangle < 0, k \in I(z_0)\}, I(z_0) = \{k \in I: W_k(z_0) = 0\}$$
is a local tent.

It should be noted that for a convex set $F \subset \mathbb{R}^{3n}$ the cone defined as follows is a local tent and exists always [21]:
$$K_F(z_0) = \{\overline{z}: \overline{z} = \lambda(z_1 - z), \lambda > 0, \forall z_1 \in F\}.$$

Section 5 concerns the following problem, labelled as (PC) for a SSDFI:

$$\text{minimize } J[x(\cdot)] = \int_0^1 f(x(t),t)dt + q(x(1)) \qquad (4)$$

(PC) $\quad W_k(x(t), x'(t), x''(t)) \leq 0$, a.e. $t \in [0,1], k \in I = \{1,...,m\}$ (5)

$$x(0) = \upsilon_0, \quad x'(0) = \upsilon_1. \qquad (6)$$

Here, firstly the convex problem is considered for simplicity; $f$ is continuously differentiable and convex with respect to $x$, $W_k, k \in I$ and $q$ are continuous convex functions and $\upsilon_0, \upsilon_1$ are fixed vectors (Remark that Theorem 5.3 extends the result for the convex problem to a more general nonconvex setting of Problem (PC)). It is required to find an arc $\tilde{x}(t)$ of the Cauchy problem (4) – (6) for the system of SSDFI satisfying the system of differential inequalities (5) almost everywhere (a.e.) on $[0,1]$ and the initial conditions (6) minimizing the Bolza functional $J[x(\cdot)]$. For this purpose, we first derive necessary optimality conditions for the discrete-approximation problem and then, by passing to the limit procedure establish sufficient optimality conditions to the problem (PC) posed by SSDFI. Here, a feasible trajectory $x(\cdot)$ is an absolutely continuous twice differentiable function on $[0,1]$ for which $x''(\cdot) \in L_1^n([0,1])$. Notice that such class of functions constitutes a Banach space with various equivalent norms.

## 3. Optimality Conditions for a SSDSI

Fisrtly, we consider the convex problem (PD). After transformation (PD) to an equivalent form by using a convex programming method and subdifferential calculus we derive an optimality condition for it.

Consider the following convex sets in the space $\mathbb{R}^{n(N+1)}$:
$$M_t = \{w = (x_0,...,x_N): W_k(x_t, x_{t+1}, x_{t+2}) \leq 0, k \in I\}, \ t = 0,1,...,N-2,$$
$$G_0 = \{w = (x_0,...,x_N): x_0 = \bar{v}_0\}, \quad G_1 = \{w = (x_0,...,x_N): x_1 = \bar{v}_1\},$$
where $w = (x_0, x_1,...,x_N) \in \mathbb{R}^{n(N+1)}$.

Denoting $g(w) = \sum_{t=2}^{N-1} f(x_t, t)$ we can replace the convex problem (PD) by the equivalent problem with geometric constraints:
$$\text{minimize } g(w) \text{ subject to } Q = G_0 \cap G_1 \cap \left(\bigcap_{t=0}^{N-2} M_t\right). \tag{7}$$

Let us compute the dual cone $K^*_{M_t}(w)$ generated by subdifferential set $\partial_z W_k(x_t, x_{t+1}, x_{t+2})$, where $K_{M_t}(w)$ is a cone of tangent directions.

**Proposition 3.1** Let $\operatorname{cone} \partial_z W_k(x_t, x_{t+1}, x_{t+2})$, $z = (x, v_1, v_2)$ be the cone generated by the subdifferential set $\partial_z W_k(x_t, x_{t+1}, x_{t+2})$. Then under the non-degeneracy condition we have
$$K^*_{M_t}(w) = \left\{w^* = (x_0^*,...,x_N^*): (x_t^*, x_{t+1}^*, x_{t+2}^*) \in \left[-\sum_{k \in I(x_t, x_{t+1}, x_{t+2})} \operatorname{cone} \partial_z W_k(x_t, x_{t+1}, x_{t+2})\right],\right.$$
$$\left. x_k^* = 0, k \neq t, t+1, t+2\right\}, \ I(x_t, x_{t+1}, x_{t+2}) = \{k \in I: W_k(x_t, x_{t+1}, x_{t+2}) = 0\}.$$

*Proof.* Obviously, by Definition 1.1 if $w + \lambda \bar{w} \in M_t$, $t = 0,...,N-2$ for sufficiently small $\lambda > 0$, that is $(x_t + \lambda \bar{x}_t, x_{t+1} + \lambda \bar{x}_{t+1}, x_{t+2} + \lambda \bar{x}_{t+2}) \in F$, $F = \{(x, v_1, v_2): W_k(x, v_1, v_2) \leq 0, k \in I\}$ then it follows that $\bar{w} \in K_{M_t}(w)$ and so
$$K_{M_t}(w) = \{\bar{w}: (\bar{x}_t, \bar{x}_{t+1}, \bar{x}_{t+2}) \in K_F(x_t, x_{t+1}, x_{t+2})\},$$
$$K_F(x_t, x_{t+1}, x_{t+2}) \equiv \operatorname{cone}[F - (x_t, x_{t+1}, x_{t+2})].$$

Recall that if $I(x_t, x_{t+1}, x_{t+2}) = \varnothing$, then $K_{M_t}(x_t, x_{t+1}, x_{t+2}) = \mathbb{R}^{3n}$ and $K^*_{M_t}(x_t, x_{t+1}, x_{t+2}) = \{0\}$. Therefore, by Theorem 1.34 [21]

$$K_F^*(x_t, x_{t+1}, x_{t+2}) = \begin{cases} \{0\} \in \mathbb{R}^{3n}, & \text{if } I(x_t, x_{t+1}, x_{t+2}) = \varnothing, \\ -\sum_{k \in I(x_t, x_{t+1}, x_{t+2})} \operatorname{cone} \partial_z W_k(x_t, x_{t+1}, x_{t+2}), & \text{if } I(x_t, x_{t+1}, x_{t+2}) \neq \varnothing. \end{cases}$$

The proof immediately follows from the arbitrariness of components $x_k$, $k \neq t, t+1, t+2$ of vectors $\bar{w}$. □

By analogy we derive the following formulas:

$$K^*_{G_0}(w) = \left\{w^* = \left(x_0^*,...,x_N^*\right): x_t^* = 0, t \neq 0\right\}, \tag{8}$$
$$K^*_{G_1}(w) = \left\{w^* = \left(x_0^*,...,x_N^*\right): x_t^* = 0, t \neq 1\right\}.$$

In the sense of the cone generated by the subdifferential set $\partial_z W_k(x_t, x_{t+1}, x_{t+2})$, we can now prove the optimality conditions for the problem (PD).

**Theorem 3.1** Let $W_k, k \in I$ be continuous convex function, $f(\cdot, t)$ $(t = 2,..., N-2)$ be continuously differentiable proper convex function in $x$, and $f(\cdot, N-1)$ be continuous convex proper function. Besides, let $f(\cdot, t)$ be continuous at points $x_t$ of some feasible trajectory. Then, for optimality of the trajectory $\{\tilde{x}_t\}_{t=0}^N$ in the problem (PD), it is necessary that there exist numbers $\mu \in \{0,1\}$, $\alpha_t^k \geq 0$, $k \in I(\tilde{x}_t, \tilde{x}_{t+1}, \tilde{x}_{t+2})$, $t = 0,..., N-2$ and vectors $x_t^*, x_N^*, u_t^*$, $t = 0,..., N-1$ not all equal to zero at the same time, satisfying the system of second order adjoint discrete inclusions generated by the given inequality constraints (2)

$$\left(x_t^* - u_t^* - \mu f_x'(\tilde{x}_t, t),\ u_{t+1}^*,\ -x_{t+2}^*\right) \in \sum_{k \in I(\tilde{x}_t, \tilde{x}_{t+1}, \tilde{x}_{t+2})} \alpha_t^k \partial_{(x,v_1,v_2)} W_k(\tilde{x}_t, \tilde{x}_{t+1}, \tilde{x}_{t+2}),\ t = 2,..., N-2,$$

$$\alpha_t^k W_k(\tilde{x}_t, \tilde{x}_{t+1}, \tilde{x}_{t+2}) = 0,\ t = 2,..., N-2,\ \alpha_t^k \geq 0,\ u_0^* = 0,\ f_x'(\tilde{x}_0, 0) = f_x'(\tilde{x}_1, 1) = 0$$

and boundary conditions
$$x_{N-1}^* - u_{N-1}^* \in \mu \partial f(\tilde{x}_{N-1}, N-1),\ x_N^* = 0.$$

If moreover the non-degeneracy condition is fulfilled, these conditions suffice for the optimality of the trajectory of $\{\tilde{x}_t\}_{t=0}^N$.

*Proof.* Clearly, if $\{\tilde{x}_t\}_{t=0}^N$ is an optimal trajectory we claim that $\tilde{w} = (\tilde{x}_0,..., \tilde{x}_N)$ is a solution of the convex mathematical programming problem (7). Then it is known (see [21]) that there exist not all zero vectors $w^*(t) \in K_{M_t}^*(\tilde{w})$, $t = 0, 1,..., N-2$, $w_0^* \in K_{G_0}^*(\tilde{w})$, $w_1^* \in K_{G_1}^*(\tilde{w})$, and a number $\mu \in \{0,1\}$, such that

$$\mu w^{0*} = \sum_{t=0}^{N-2} w^*(t) + w_0^* + w_1^*,\quad w^{0*} \in \partial g(\tilde{w}). \tag{9}$$

From definition of the function $g$ we claim that the vector $w^{0*}$ has a form $w^{0*} = \left(x_{00}^*, x_{10}^*,..., x_{N0}^*\right)$, $x_{t0}^* = f_x'(\tilde{x}_t, t)$ $(t = 0,..., N-2)$, $x_{N-1,0}^* \in \partial_x f(\tilde{x}_{N-1}, N-1)$. By Proposition 3.1 and formulas (8), it is deduced that

$$w^*(t) = \left(0,...,0, x_t^*(t), x_{t+1}^*(t), x_{t+2}^*(t), 0,..., 0\right),$$
$$\left(-x_t^*(t), -x_{t+1}^*(t), -x_{t+2}^*(t)\right) \in \sum_{k \in I(\tilde{x}_t, \tilde{x}_{t+1}, \tilde{x}_{t+2})} \text{cone}\ \partial_{(x,v_1,v_2)} W_k(\tilde{x}_t, \tilde{x}_{t+1}, \tilde{x}_{t+2}),\ t = 2,..., N-2, \tag{10}$$

$I(\tilde{x}_t, \tilde{x}_{t+1}, \tilde{x}_{t+2}) = \{k \in I: W_k(\tilde{x}_t, \tilde{x}_{t+1}, \tilde{x}_{t+2}) = 0\}$, $w_0^* = \left(\overline{x}^*, 0,..., 0\right)$, $w_1^* = \left(0, \overline{\overline{x}}^*, 0,..., 0\right)$,

where $\overline{x}^*, \overline{\overline{x}}^*$ are arbitrary vectors. Besides, by the component-wise representation of (9), we have

$$\mu x_{t0}^* = x_t^*(t) + x_t^*(t-1) + x_t^*(t-2),\ t = 2,..., N-2,$$
$$\overline{x}^* + x_0^*(0) = 0,\quad \overline{\overline{x}}^* + x_1^*(1) + x_1^*(0) = 0. \tag{11}$$

Recalling that

$$\text{cone}\,\partial_{(x,v_1,v_2)}W_k(x_t, x_{t+1}, x_{t+2}) = \{(\alpha_t^k \hat{x}_t^*, \alpha_t^k \hat{x}_{t+1}^*, \alpha_t^k \hat{x}_{t+2}^*) :$$
$$(\hat{x}_t^*, \hat{x}_{t+1}^*, \hat{x}_{t+2}^*) \in \partial_{(x,v_1,v_2)}W_k(x_t, x_{t+1}, x_{t+2}),\ \alpha_t^k W_k(x_t, x_{t+1}, x_{t+2}) = 0,\ \alpha_t^k \geq 0\}$$

by the second formula of (10) we conclude that

$$-x_t^*(t) = \alpha_t^k \hat{x}_t^*,\quad -x_{t+1}^*(t) = \alpha_t^k \hat{x}_{t+1}^*,\quad -x_{t+2}^*(t) = \alpha_t^k \hat{x}_{t+2}^*,\ t=2,\ldots,N-2,\ \alpha_t^k \geq 0,$$
$$(\hat{x}_t^*, \hat{x}_{t+1}^*, \hat{x}_{t+2}^*) \in \sum_{k \in I(\tilde{x}_t, \tilde{x}_{t+1}, \tilde{x}_{t+2})} \partial_{(x,v_1,v_2)}W_k(\tilde{x}_t, \tilde{x}_{t+1}, \tilde{x}_{t+2}),\quad \alpha_t^k W_k(\tilde{x}_t, \tilde{x}_{t+1}, \tilde{x}_{t+2}) = 0. \quad (12)$$

Since $\alpha_t^k \hat{x}_t^* = -x_t^*(t)$, $\alpha_t^k \hat{x}_{t+1}^* = -x_{t+1}^*(t)$, $\alpha_t^k \hat{x}_{t+2}^* = -x_{t+2}^*(t)$ in turn from the second formula of (12) it follows that

$$(-x_t^*(t), -x_{t+1}^*(t), -x_{t+2}^*(t)) \in \sum_{k \in I(\tilde{x}_t, \tilde{x}_{t+1}, \tilde{x}_{t+2})} \alpha_t^k \partial_{(x,v_1,v_2)}W_k(\tilde{x}_t, \tilde{x}_{t+1}, \tilde{x}_{t+2}),$$
$$\alpha_t^k W_k(\tilde{x}_t, \tilde{x}_{t+1}, \tilde{x}_{t+2}) = 0,\ \alpha_t \geq 0,\ t=2,\ldots,N-2. \quad (13)$$

By employing following notations $x_{t+1}^*(t) \equiv -u_{t+1}^*$, $x_{t+2}^*(t) \equiv x_{t+2}^*$, $t=1,\ldots,N-2$ in the first formula of (11), we derive from (13) that

$$(x_t^* - u_t^* - \mu x_{t0}^*,\ u_{t+1}^*,\ -x_{t+2}^*) \in \sum_{k \in I(\tilde{x}_t, \tilde{x}_{t+1}, \tilde{x}_{t+2})} \alpha_t^k \partial_{(x,v_1,v_2)}W_k(\tilde{x}_t, \tilde{x}_{t+1}, \tilde{x}_{t+2}),\ t=2,\ldots,N-2,$$
$$\alpha_t W_k(\tilde{x}_t, \tilde{x}_{t+1}, \tilde{x}_{t+2}) = 0,\ \alpha_t \geq 0,\ t=2,\ldots,N-2,\ k \in I(\tilde{x}_t, \tilde{x}_{t+1}, \tilde{x}_{t+2}). \quad (14)$$

Moreover, by setting $f(x_0, 0) = f(x_1, 1) = 0$, $\bar{x}^* = -x_0^*$, $u_0^* \equiv 0$ and $\bar{\bar{x}}^* = -x_1^*$ in the second and third equalities, respectively the formula (14) can be generalized to the case $t=0,1$. Consequently for $t = N-1$, we get

$$\mu x_{(N-1)0}^* = x_{N-1}^* - u_{N-1}^*. \quad (15)$$

Note that since $f(\tilde{x}_N, N) \equiv 0$ it follows that $x_N^* = 0$.

Hence, considering the formulas (14) and (15), we deduce that the indicated conditions of theorem are necessary for optimality.

On the other hand, by Theorem 3.3 [21], under the nondegeneracy condition, the equality (9) holds when $\mu = 1$, which completes the sufficiency. □

Below, we show how the result of Theorem 3.1 can be extended to the problem (PD) in the nonconvex case.

**Theorem 3.2** Suppose that the condition $H_2$ for the nonconvex problem (PD) is satisfied; $W_k, k \in I$ and $f(\cdot, N-1)$ admit a continuous CUAs $h_{W_k}^t(\cdot, \tilde{x}_t, \tilde{x}_{t+1}, \tilde{x}_{t+2})$, $h_{N-1}(\cdot, \tilde{x}_{N-1})$, respectively, that is $\partial_{(x,v_1,v_2)}W_k(\tilde{x}_t, \tilde{x}_{t+1}, \tilde{x}_{t+2}) = \partial h_{W_k}^t(0, \tilde{x}_t, \tilde{x}_{t+1}, \tilde{x}_{t+2})$, $\partial f(\tilde{x}_{N-1}, N-1)$, $= \partial h_{N-1}(0, \tilde{x}_{N-1})$ and $f(\cdot, t)$ $(t = 2,\ldots,N-2)$ is continuously differentiable nonconvex function. In addition, let $\bar{z}$ be a point such that $h_{W_k}^t(\bar{z}, \tilde{x}_t, \tilde{x}_{t+1}, \tilde{x}_{t+2}) < 0$ and $W_k(\tilde{x}_t, \tilde{x}_{t+1}, \tilde{x}_{t+2}) = 0$. Then, for optimality of the trajectory $\{\tilde{x}_t\}_{t=0}^N$ in the nonconvex problem with the SSDSI, it is necessary that there exist numbers $\mu \geq 0$, $\alpha_t^k \geq 0, k \in I(\tilde{x}_t, \tilde{x}_{t+1}, \tilde{x}_{t+2})$, $t = 2,...,N-2$ and a pair of vectors $\{x_t^*\}, \{u_t^*\}$, not all equal to zero, fulfilling the analogous conditions in the previous case.

*Proof.* In the nonconvex problem (PD), the assumption $H_2$ establishes the conditions of Theorem 3.24 [21] for the equivalent problem (7). So, by analogy with this theorem for the main relation (9), the necessary condition is found as in Theorem 3.1. Note that in this case the required cone

$$K_F(\tilde{x}_t, \tilde{x}_{t+1}, \tilde{x}_{t+2}) = \{\bar{z} : h_{W_k}^t(\bar{z}, \tilde{x}_t, \tilde{x}_{t+1}, \tilde{x}_{t+2}) < 0, k \in I(\tilde{x}_t, \tilde{x}_{t+1}, \tilde{x}_{t+2})\}$$

is a local tent [21]. □

**Remark 3.1** Because of the arbitrariness of $\bar{x}^*, \bar{\bar{x}}^*$ ($\bar{x}^* = -x_0^*, \bar{\bar{x}}^* = -x_1^*$) it is not hard to see that the second and third equalities of formulas (11) are realized always. Then, it is clear that the conditions $u_0^* = 0$, $f_x'(\tilde{x}_0, 0) = f_x'(\tilde{x}_1, 1) = 0$ can be disregarded without loss of generality and second order discrete adjoint inclusions are preserved for $t = 2, 3,..., N-2$.

**Remark 3.2** Another approach to solving the problem (PD) is to represent the system of inequalities (2) using a set-valued mapping

$$F(x, v_1) = \{v_2 : W_k(x, v_1, v_2) \le 0, k \in I\}.$$

Then the condition (2) can be reduced to the following discrete inclusions

$$x_{t+2} \in F(x_t, x_{t+1}), t = 0,..., N-2.$$

Now, suppose that $W_k$, $k \in I$ are convex functions, continuous at $(x^0, v_1^0)$, where $W_k(x^0, v_1^0) = 0, k \in I_0 \equiv I(x^0, v_1^0) = \{k : W_k(x^0, v_1^0) = 0\}$. Besides, assume that there is a point $(x^1, v_1^1)$ such that $W_k(x^1, v_1^1) < 0$. Then by Theorem 2.13 [21] the locally adjoint mapping has the form

$$F^*(v_2^*;(x^0, v_1^0)) = \left\{-\sum_{k \in I_0} \alpha^k x_k^* : v_1^* = \sum_{k \in I_0} \alpha^k v_1^{k*}, (x_k^*, v_1^{k*}) \in \partial_{(x, v_1)} W_k(x^0, v_1^0), \alpha^k \ge 0, k \in I_0\right\}.$$

Then, using a locally conjugate multi-valued mapping apparatus $F^*$ (see, for example, [27-29]), we can derive necessary and sufficient optimality conditions for Theorem 3.1.

## 4. Equivalence of Subdifferentials and Optimization of Discrete-Approximation Problem with SSDSI

Recall the Theorem 3.1 [22] important for what follows. By this theorem for proper convex function $\Phi_k(x, v_1, v_2) \equiv W_k\left[x, \frac{1}{\delta}(v_1 - x), \frac{1}{\delta^2}(v_2 - 2v_1 + x)\right]$, the following inclusions are equivalent:

(1) $(\bar{x}^*, \bar{v}_1^*, \bar{v}_2^*) \in \partial_{(x, v_1, v_2)} \Phi_k(z_0), z_0 = (x^0, v_1^0, v_2^0) \in \text{dom}\, \Phi_k, k \in I,$

(2) $(\bar{x}^* + \bar{v}_1^* + \bar{v}_2^*, \delta\bar{v}_1^* + 2\delta\bar{v}_2^*, \delta^2\bar{v}_2^*) \in \partial W_k\left[x^0, \frac{1}{\delta}(v_1^0 - x^0), \frac{1}{\delta^2}(v_2^0 - 2v_1^0 + x^0)\right].$ (16)

**Proposition 4.1** For a proper convex function given by $\Phi_k(x, v_1, v_2) \equiv W_k^1\left[x, \frac{1}{\delta^2}(v_2 - 2v_1 + x)\right]$, the following inclusions under the condition that $\bar{v}_1^* = -2\bar{v}_2^*$ are equivalent:

(1) $(\bar{x}^*, \bar{v}_1^*, \bar{v}_2^*) \in \partial_z \Phi_k(z_0)$, $z_0 = (x^0, v_1^0, v_2^0) \in \mathrm{dom}\Phi_k, k \in I$,

(2) $(\bar{x}^* - \bar{v}_2^*, \delta^2 \bar{v}_2^*) \in \partial W_k^1\left[x^0, \frac{1}{\delta^2}(v_2^0 - 2v_1^0 + x^0)\right]$.

*Proof.* Note that $\partial_z \Phi_k(z_0)$ is a convex closed set and is bounded for $z_0 \in \mathrm{ri}(\mathrm{dom}\Phi_k)$ [21]. Furthermore, by definition of subdifferential one has

$$\partial_z \Phi_k(z_0) = \{(\bar{x}^*, \bar{v}_1^*, \bar{v}_2^*): \Phi_k(z) - \Phi_k(z_0) \geq \langle \bar{x}^*, x - x^0 \rangle \\ + \langle \bar{v}_1^*, v_1 - v_1^0 \rangle + \langle \bar{v}_2^*, v_2 - v_2^0 \rangle, \forall z = (x, v_1, v_2) \in \mathbb{R}^{3n}\} \tag{17}$$

and similarly

$$\partial W_k^1\left[x^0, \frac{1}{\delta^2}(v_2^0 - 2v_1^0 + x^0)\right] = \left\{(x^*, v^*): W_k^1\left[x, \frac{1}{\delta^2}(v_2 - 2v_1 + x)\right] \right.$$
$$\left. - W_k^1\left[x^0, \frac{1}{\delta^2}(v_2^0 - 2v_1^0 + x^0)\right] \geq \langle x^*, x - x^0 \rangle + \left\langle v^*, \frac{1}{\delta^2}(v_2 - 2v_1 + x) - \frac{1}{\delta^2}(v_2^0 - 2v_1^0 + x^0) \right\rangle, \right.$$
$$\forall z = (x, v_1, v_2) \in \mathbb{R}^{3n}.$$

After obvious simplification we have

$$\partial W_k^1\left[x^0, \frac{1}{\delta^2}(v_2^0 - 2v_1^0 + x^0)\right] = \left\{(x^*, v^*): \left\langle x^* + \frac{v^*}{\delta^2}, x - x^0 \right\rangle \right.$$
$$\left. + \left\langle -\frac{2v^*}{\delta^2}, v_1 - v_1^0 \right\rangle + \left\langle \frac{v^*}{\delta^2}, v_2 - v_2^0 \right\rangle, \forall (x, v_1, v_2) \in \mathbb{R}^{3n} \right\}. \tag{18}$$

By comparing (17) and (18), we can derive that

$$\bar{x}^* = x^* + \frac{v^*}{\delta^2}; \quad \bar{v}_1^* = -\frac{2v^*}{\delta^2}; \quad \bar{v}_2^* = \frac{v^*}{\delta^2},$$

whence

$$x^* = \bar{x}^* - \bar{v}_2^*, \quad \bar{v}_1^* = -2\bar{v}_2^*, \quad v_2^* = \delta^2 \bar{v}_2^*.$$

This means that (1) and (2) of theorem are equivalent. □

Following is an alteration of Theorem 3.1 [22].

**Proposition 4.2** For a proper convex function given by $\Phi_k(x, v_1, v_2) \equiv W_k^2\left[x, \frac{1}{\delta}(v_1 - x)\right]$ under the condition that $\bar{v}_2^* \equiv 0$, the following inclusions are equivalent:

$(1^0)$ $(\bar{x}^*, \bar{v}_1^*, \bar{v}_2^*) \in \partial_z \Phi_k(z_0)$, $z_0 = (x^0, v_1^0, v_2^0) \in \mathrm{dom}\Phi_k$,

$(2^0)$ $(\bar{x}^* + \bar{v}_1^*, \delta \bar{v}_1^*) \in \partial W_k^2\left[x^0, \frac{1}{\delta}(v_1^0 - x^0)\right]$.

*Proof.* Clearly, for all $(x, v_1) \in \mathbb{R}^{2n}$ we have

$$\partial W_k^2\left[x^0, \frac{1}{\delta}(v_1^0 - x^0)\right] = \left\{(x^*, v^*) : W_k^2\left[x, \frac{1}{\delta}(v_1 - x)\right] - W_k^2\left[x^0, \frac{1}{\delta}(v_1^0 - x^0)\right]\right.$$

$$\left. \geq \left\langle x^* - \frac{v^*}{\delta}, x - x^0 \right\rangle + \left\langle \frac{v^*}{\delta}, v_1 - v_1^0 \right\rangle \right\}. \tag{19}$$

Then, the proof of proposition follows immediately from relations (17) and (19). □

In the next parts, it will be explained how to construct the discrete-approximation problem for (PC); but first consider the first and second order difference operators:

$$\Delta x(t) = \frac{1}{\delta}[x(t+\delta) - x(t)], \quad \Delta^2 x(t) = \frac{1}{\delta}[\Delta x(t+\delta) - \Delta x(t)], \quad t = 0, \delta, \ldots, 1 - 2\delta,$$

where $\delta$ is a step on the $t$-axis and $x(t)$ is a grid function on a uniform grid on $[0,1]$.

Let us define following discrete-approximation problem associated with the problem (PC):

$$\text{minimize } J_\delta[x(\cdot)] = \sum_{t=2\delta,\ldots,1-2\delta} \delta f(x(t),t) + q(x(1-\delta)),$$

subject to

$$W_k\left(x(t), \Delta x(t), \Delta^2 x(t)\right) \leq 0, \quad t = 0, \delta, 2\delta, \ldots, 1-2\delta, \quad k \in I, \tag{20}$$

$$x(0) = v_0, \quad \Delta x(0) = v_1.$$

To utilize the results of Theorem 3.1 in the problem (20), we reduce the problem (20) to a problem of the form (PDA):

$$\text{minimize } J_\delta[x(\cdot)] = \sum_{t=2\delta,\ldots,1-2\delta} \delta f(x(t),t) + q(x(1-\delta)), \tag{21}$$

(PDA)

$$\Phi_k\left(x(t), x(t+\delta), x(t+2\delta)\right) \leq 0, \quad k \in I, \quad t = 0, \delta, 2\delta, \ldots, 1-2\delta, \tag{22}$$

$$x(0) = v_0, \quad x(\delta) = v_0 + \delta v_1,$$

where

$$\Phi_k\left(x(t), x(t+\delta), x(t+2\delta)\right) \equiv W_k\left(x(t), \Delta x(t), \Delta^2 x(t)\right). \tag{23}$$

An immediate useful consequence of the equivalence theorems, key instruments in the study of discrete-approximation problems, is the following result.

**Theorem 4.1** Let $W_k$ $(k \in I)$, $q$ be continuous convex functions, and $f(\cdot, t)$ $(t = 2\delta, 3\delta, \ldots, 1 - 2\delta)$ be continuously differentiable proper functions convex with respect to $x$. Besides, $f(\cdot, t)$ is continuous at points $x_t$ of some feasible trajectory. Then, for optimality of the trajectory $\{\tilde{x}(t)\}$ in the discrete-approximate problem (21), (22) it is necessary that there exist numbers $\mu = \mu_\delta \in \{0,1\}$, $\alpha_k(t) \geq 0$, $k \in I(\tilde{x}(t), \Delta \tilde{x}(t), \Delta^2 \tilde{x}(t))$ and a pair of vectors $\{x^*(t), u^*(t)\}$, not all equal to zero at the same time, satisfying the second order adjoint inclusions generated by inequality constraints (22)

$$\left[\frac{1}{\delta^2}\left(x^*(t)-u^*(t)+u^*(t+\delta)-x^*(t+2\delta)-\mu\delta f_x'(\tilde{x}(t),t)\right),\right.$$

$$\left.\frac{1}{\delta}\left(u^*(t+\delta)-2x^*(t+2\delta)\right),\ -x^*(t+2\delta)\right] \in \sum_{k\in I(\tilde{x}(t),\Delta\tilde{x}(t),\Delta^2\tilde{x}(t))} \alpha_k(t)\partial_{(x,v_1,v_2)}W_k(\tilde{x}(t),\Delta\tilde{x}(t),\Delta^2\tilde{x}(t)),$$

$$\alpha_k(t)W_k(\tilde{x}(t),\Delta\tilde{x}(t),\Delta^2\tilde{x}(t)) = 0,\ \alpha_k(t)\geq 0,\ t=2\delta, 3\delta,...,1-2\delta \quad (24)$$

and boundary condition

$$x^*(1-\delta)-u^*(1-\delta) \in \mu\partial q\big(\tilde{x}(1-\delta)\big),\quad x^*(1) = 0. \quad (25)$$

In addition, under the nondegeneracy condition, these conditions suffice for optimality of $\{\tilde{x}(t)\}$.

*Proof.* The proof is immediate from the results of Theorem 3.1. Indeed by Theorem 3.1 for optimality of the trajectory $\{\tilde{x}(t)\}:=\{\tilde{x}(t):t=0,\delta,...,1\}$, in problem (21), (22) it is necessary that there exist a pair of vectors $\{u^*(t), x^*(t)\}$ and numbers $\mu=\mu_\delta\in\{0,1\}$, $\alpha_k(t)\geq 0$ not all zero, such that

$$\big(x^*(t)-u^*(t)-\mu\delta f_x'(\tilde{x}(t),t),\ u^*(t+\delta),\ -x^*(t+2\delta)\big)$$
$$\in \sum_{k\in I(\tilde{x}(t),\tilde{x}(t+\delta),\tilde{x}(t+2\delta))} \alpha_k(t)\partial_z\Phi_k\big(\tilde{x}(t),\tilde{x}(t+\delta),\tilde{x}(t+2\delta)\big) \quad (26)$$

$$\alpha_k(t)\Phi_k(\tilde{x}(t),\tilde{x}(t+\delta),\tilde{x}(t+2\delta))=0,\ \alpha_k(t)\geq 0,$$

$$t=2\delta,...,1-2\delta,\ u^*(0)=0,\ f_x'(\tilde{x}(0),0)=f_x'(\tilde{x}(\delta),\delta)=0,$$

$$x^*(1-\delta)-u^*(1-\delta)\in\mu\partial q(\tilde{x}(1-\delta)),\ x^*(1)=0. \quad (27)$$

So it remains only to express in the relationship (26) the subdifferential $\partial_z\Phi_k(\tilde{x}(t),\tilde{x}(t+\delta),\tilde{x}(t+2\delta))$ in terms of subdifferential $\partial_z W_k(\tilde{x}(t),\Delta\tilde{x}(t),\Delta^2\tilde{x}(t))$. Indeed by Theorem 3.1 and formula (16), the second order adjoint inclusions (26) for convex problem take the form (24), where the values $\alpha(t)/\delta^2$, $\mu/\delta$ are shown again by $\alpha(t)$ and $\mu$, respectively. □

**Remark 4.1** Similarly, using Theorem 3.2 we can generalize the result of Theorem 4.1 to nonconvex problem (PDA). Assume that the condition $H_2$ is satisfied for the nonconvex problem (PDA). Then, for optimality of the trajectory $\{\tilde{x}(t)\}$ it is necessary that there exist numbers $\mu\geq 0$, $\alpha_k(t)\geq 0$, $k\in I(\tilde{x}(t),\Delta\tilde{x}(t),\Delta^2\tilde{x}(t))$ and a pair of vectors $\{x^*(t),u^*(t)\}$, not all equal to zero at the same time, satisfying the system of second order adjoint inclusions (24) and boundary condition (25) rewritten for nonconvex case.

Later on we shall investigate the discrete problems for the different particular case of $W_k$ function.

**Theorem 4.2** If in the problem (PC), $W_k(x, x', x'') \equiv W_k^1(x, x'')$, i.e. a function $W_k$ does not depend on the second argument $x'$, then for the corresponding (PDA) problem, the conditions (24), (25) have the forms

$$\left[\Delta^2 x^*(t) - \mu f_x'(\tilde{x}(t), t), -x^*(t+2\delta)\right] \in \sum_{k \in I(\tilde{x}(t), \Delta^2 \tilde{x}(t))} \alpha_k(t) \partial_{(x, v_2)} W_k^1(\tilde{x}(t), \Delta^2 \tilde{x}(t))$$

$$x^*(1-\delta) \in \mu \partial q(\tilde{x}(1-\delta)), \quad x^*(1) = 0, \quad t = 2\delta, 3\delta, \ldots, 1-2\delta,$$

respectively.

*Proof.* Applying the Proposition 4.1 in the second order adjoint inclusions (26) we have

$$\left[\frac{1}{\delta^2}\left(x^*(t) - u^*(t) + u^*(t+\delta) - x^*(t+2\delta) - \mu \delta f_x'(\tilde{x}(t), t)\right), \right.$$
$$\left. -x^*(t+2\delta)\right] \in \sum_{k \in I(\tilde{x}(t), \Delta^2 \tilde{x}(t))} \alpha_k(t) \partial_{(x, v_2)} W_k(\tilde{x}(t), \Delta^2 \tilde{x}(t)). \tag{28}$$

On the other hand again by Proposition 4.1 $\bar{v}_1^* = -2\bar{v}_2^*$, which by second order adjoint inclusions (26) means that $u^*(t) = 2x^*(t+\delta)$. Then taking into account the last condition in the left hand side of (28) we have

$$\frac{1}{\delta^2}\left[x^*(t) - u^*(t) + u^*(t+\delta) - x^*(t+2\delta)\right] = \Delta^2 x^*(t).$$

Finally, by substitution $u^*(1-\delta) = 2x^*(1)$ into (27), we get the desirable result. □

**Theorem 4.3** Suppose that in the problem (PC) $W_k(x, x', x'') \equiv W_k^2(x, x')$, i.e. a function $W_k$ doesn't depends on the third argument $x''$. Then for the corresponding (PDA) problem the second order adjoint inclusions (24) and boundary condition (25) have the forms

$$\left[\left(\Delta u^*(t) - \mu f_x'(\tilde{x}(t), t), u^*(t+\delta)\right)\right] \in \sum_{k \in I(\tilde{x}(t), \Delta \tilde{x}(t))} \alpha_k(t) \partial_{(x, v_1)} W_k^2(\tilde{x}(t), \Delta \tilde{x}(t)), \tag{29}$$

$$-u^*(1-\delta) \in \mu \partial q(\tilde{x}(1-\delta)), \quad x^*(1) = 0, \quad t = 2\delta, 3\delta, \ldots, 1-2\delta.$$

*Proof.* Using the Proposition 4.2 and inclusion (26) and considering that $x^*(1-\delta) = 0$ ($\bar{v}_2^* \equiv 0$) we have

$$\left[\frac{1}{\delta}\left(u^*(t+\delta) - u^*(t) - \mu \delta f_x'(\tilde{x}(t), t)\right), u^*(t+\delta)\right] \in \sum_{k \in I(\tilde{x}(t), \Delta \tilde{x}(t))} \alpha_k(t) \partial_{(x, v_1)} W_k^2(\tilde{x}(t), \Delta \tilde{x}(t)), \tag{30}$$

which implies the inclusion (29). On the other hand, by substitution $x^*(t) \equiv 0$ into (25) we have $-u^*(1-\delta) \in \mu \partial q(\tilde{x}(1-\delta))$. □

Let us now rewrite (24), (25) in more convenient form.

**Lemma 4.1** The second order adjoint inclusions (24) and boundary condition (25) have the forms

$$\left[\Delta^2 x^*(t) + \psi^*(t) - \mu f_x'(\tilde{x}(t),t),\ \psi^*(t),\ -x^*(t+2\delta)\right]$$
$$\in \sum_{k \in I(\tilde{x}(t),\Delta\tilde{x}(t),\Delta^2\tilde{x}(t))} \alpha_k(t)\partial_{(x,v_1,v_2)} W_k(\tilde{x}(t),\Delta\tilde{x}(t),\Delta^2\tilde{x}(t)),$$
$$\psi^*(1) + \Delta x^*(1-\delta) \in \mu\partial q(\tilde{x}(1-\delta)),\ x^*(1) = 0,\ t = 2\delta, 3\delta, \ldots, 1-2\delta,$$
$$\psi^*(t) = \frac{1}{\delta}\left[u^*(t) - 2x^*(t+\delta)\right],$$

respectively.

*Proof.* The condition $u^*(t+\delta) = \delta\psi^*(t+\delta) + 2x^*(t+2\delta)$ yields immediately

$$\frac{1}{\delta^2}\left[x^*(t) - u^*(t) + u^*(t+\delta) - x^*(t+2\delta)\right] = \Delta^2 x^*(t) + \Delta\psi^*(t). \tag{31}$$

Moreover, by the boundary condition (25) one has

$$\frac{1}{\delta}\left[x^*(1-\delta) + 2x^*(1) - x^*(1) - u^*(1-\delta)\right] \in \mu\partial q(\tilde{x}(1-\delta))$$

or consequently

$$-\psi^*(1-\delta) - \Delta x^*(1-\delta) \in \mu\partial q(\tilde{x}(1-\delta)). \tag{32}$$

Thus, taking into account (31), (32) in the second order adjoint inclusions (24) and boundary condition (25) respectively, we get the desirable result. □

## 5. Sufficient Conditions of Optimality for SSDFI

By means of the results in Section 4, we can assert a sufficient condition for optimality of the continuous problem (PC). In fact, setting $\mu = 1$ in conditions of Lemma 4.1 and as $\delta \to 0$, we establish the second order adjoint differential inclusions generated by the given SSDFI (5):

(i)
$$\left(\frac{d^2 x^*(t)}{dt^2} + \frac{d\psi^*(t)}{dt} - f_x'(\tilde{x}(t),t),\ \psi^*(t),\ -x^*(t)\right)$$
$$\in \sum_{k \in I(\tilde{x}(t),\tilde{x}'(t),\tilde{x}''(t))} \alpha_k(t)\partial_{(x,v_1,v_2)} W_k(\tilde{x}(t),\tilde{x}'(t),\tilde{x}''(t)),\ \text{a.e.}\ t \in [0,1],$$
$$I(\tilde{x}(t),\tilde{x}'(t),\tilde{x}''(t)) = \{k \in I:\ W_k(\tilde{x}(t),\tilde{x}'(t),\tilde{x}''(t)) = 0\}.$$

Along with this we have

(ii) $\quad \alpha_k(t)W_k\left(\tilde{x}(t),\tilde{x}'(t),\tilde{x}''(t)\right) = 0,\ \alpha_k(t) \geq 0$, a.e. $t \in [0,1]$

and the boundary condition

(iii) $\quad -\psi^*(1) - \dfrac{dx^*(1)}{dt} \in \partial q(\tilde{x}(1)),\ x^*(1) = 0.$

Next, suppose that $x^*(t)$, $t \in [0,1]$ is an absolutely continuous, twice differentiable function on $[0,1]$ for which $x^{*\prime\prime}(\cdot) \in L_1^n([0,1])$. Moreover suppose that $\psi^*(t)$, $t \in [0,1]$ is an absolutely continuous and $\psi^{*\prime}(\cdot) \in L_1^n([0,1])$.

The next theorem suggests a convenient way to check whether a pair of absolutely continuous functions $\{x^*(t), \psi^*(t)\}$ is optimal or not.

**Theorem 5.1** Let $f(\cdot, t): \mathbb{R}^n \to \mathbb{R}^1$ be continuously differentiable convex with respect to $x$ function, and $W_k (k \in I)$, $q$ be continuous convex functions. Then, for optimality of the arc $\tilde{x}(t)$ in the convex problem (PC) it is sufficient that there exist a pair of absolutely continuous functions $\{x^*(t), \psi^*(t)\}$, $t \in [0,1]$ and function $\alpha_k(t) \geq 0$, $k \in I(\tilde{x}(t), \tilde{x}'(t), \tilde{x}''(t))$ satisfying a.e. the second order adjoint differential inclusion generated by given SSDFI (i), (ii) and boundary condition (iii).

*Proof.* By the definition of the subdifferential of the function $W_k$, we can rewrite the second order adjoint differential inclusion (i) as:

$$\sum_{k \in I(\tilde{x}(t), \tilde{x}'(t), \tilde{x}''(t))} \alpha_k(t) W_k(x(t), x'(t), x''(t)) - \sum_{k \in I(\tilde{x}(t), \tilde{x}'(t), \tilde{x}''(t))} \alpha_k(t) W_k(\tilde{x}(t), \tilde{x}'(t), \tilde{x}''(t)) \quad (33)$$

$$\geq \left\langle \frac{d^2 x^*(t)}{dt^2} + \frac{d\psi^*(t)}{dt} - f'_x(\tilde{x}(t), t), x(t) - \tilde{x}(t) \right\rangle$$

$$+ \left\langle \psi^*(t), \frac{d(x(t) - \tilde{x}(t))}{dt} \right\rangle - \left\langle x^*(t), \frac{d^2(x(t) - \tilde{x}(t))}{dt^2} \right\rangle.$$

Besides, since $f(\cdot, t)$ is convex for all feasible solutions $x(\cdot)$

$$\langle -f'_x(\tilde{x}(t), t), x(t) - \tilde{x}(t) \rangle \geq f(\tilde{x}(t), t) - f(x(t), t)$$

the inequality (33) can be reduced to

$$\sum_{k \in I(\tilde{x}(t), \tilde{x}'(t), \tilde{x}''(t))} \left[ \alpha_k(t) W_k(x(t), x'(t), x''(t)) - \alpha_k(t) W_k(\tilde{x}(t), \tilde{x}'(t), \tilde{x}''(t)) \right]$$

$$\geq \left\langle \frac{d^2 x^*(t)}{dt^2} + \frac{d\psi^*(t)}{dt}, x(t) - \tilde{x}(t) \right\rangle + f(\tilde{x}(t), t) - f(x(t), t)$$

$$+ \left\langle \psi^*(t), \frac{d(x(t) - \tilde{x}(t))}{dt} \right\rangle - \left\langle x^*(t), \frac{d^2(x(t) - \tilde{x}(t))}{dt^2} \right\rangle.$$

Now because by the condition (ii) $\alpha_k(t) W_k(\tilde{x}(t), \tilde{x}'(t), \tilde{x}''(t)) = 0$ ($\alpha_k(t) \geq 0$) and $\alpha_k(t) W_k(x(t), x'(t), x''(t)) \leq 0$ for all feasible solutions $x(\cdot)$, it follows from the last inequality

$$0 \geq \left\langle \frac{d^2 x^*(t)}{dt^2} + \frac{d\psi^*(t)}{dt}, x(t) - \tilde{x}(t) \right\rangle + f(\tilde{x}(t),t) - f(x(t),t)$$
$$+ \left\langle \psi^*(t), \frac{d(x(t) - \tilde{x}(t))}{dt} \right\rangle - \left\langle x^*(t), \frac{d^2(x(t) - \tilde{x}(t))}{dt^2} \right\rangle.$$

For ease, let us rewrite this inequality as

$$f(x(t),t) - f(\tilde{x}(t),t) \geq \left\langle \frac{d^2 x^*(t)}{dt^2}, x(t) - \tilde{x}(t) \right\rangle \tag{34}$$

$$- \left\langle \frac{d^2(x(t) - \tilde{x}(t))}{dt^2}, x^*(t) \right\rangle + \frac{d}{dt} \left\langle \psi^*(t), x(t) - \tilde{x}(t) \right\rangle.$$

Now considering that $x(\cdot), \tilde{x}(\cdot)$ are feasible ($x(0) = \tilde{x}(0) = \upsilon_0$) by integrating the inequality (34) over $[0,1]$ yields

$$\int_0^1 [f(x(t),t) - f(\tilde{x}(t),t)] dt \geq \int_0^1 \left[ \left\langle \frac{d^2 x^*(t)}{dt^2}, x(t) - \tilde{x}(t) \right\rangle - \left\langle \frac{d^2(x(t) - \tilde{x}(t))}{dt^2}, x^*(t) \right\rangle \right] dt \tag{35}$$

$$+ \left\langle \psi^*(1), x(1) - \tilde{x}(1) \right\rangle - \left\langle \psi^*(0), x(0) - \tilde{x}(0) \right\rangle$$

$$= \int_0^1 \left[ \left\langle \frac{d^2 x^*(t)}{dt^2}, x(t) - \tilde{x}(t) \right\rangle - \left\langle \frac{d^2(x(t) - \tilde{x}(t))}{dt^2}, x^*(t) \right\rangle \right] dt + \left\langle \psi^*(1), x(1) - \tilde{x}(1) \right\rangle.$$

Denoting by $\Omega$ the expression in the brackets on the right hand side of (35), we can rewrite as

$$\Omega = \frac{d}{dt} \left\langle \frac{dx^*(t)}{dt}, x(t) - \tilde{x}(t) \right\rangle - \frac{d}{dt} \left\langle \frac{d(x(t) - \tilde{x}(t))}{dt}, x^*(t) \right\rangle.$$

Then, it can be easily verified that

$$\int_0^1 \Omega dt = \left\langle \frac{dx^*(1)}{dt}, x(1) - \tilde{x}(1) \right\rangle - \left\langle \frac{dx^*(0)}{dt}, x(0) - \tilde{x}(0) \right\rangle \tag{36}$$

$$- \left\langle \frac{d(x(1) - \tilde{x}(1))}{dt}, x^*(1) \right\rangle + \left\langle \frac{d(x(0) - \tilde{x}(0))}{dt}, x^*(0) \right\rangle.$$

Because by condition (iii) $x^*(1) = 0$ and $x(t), \tilde{x}(t)$ are feasible solutions ($x(0) = \tilde{x}(0) = \upsilon_0$, $x'(0) = \tilde{x}'(0) = \upsilon_1$) the relation (36) can be expressed as

$$\int_0^1 \Omega dt = \left\langle \frac{dx^*(1)}{dt}, x(1) - \tilde{x}(1) \right\rangle. \tag{37}$$

Thus, (35) and (37) imply that

$$\int_0^1 [f(x(t),t) - f(\tilde{x}(t),t)] dt \geq \left\langle \frac{dx^*(1)}{dt}, x(1) - \tilde{x}(1) \right\rangle + \left\langle \psi^*(1), x(1) - \tilde{x}(1) \right\rangle. \tag{38}$$

Now remember that by the condition (iii) for all feasible arcs $x(t)$, $t \in [0,1]$, we have

$$q(x(1)) - q(\tilde{x}(1)) \geq -\left\langle \psi^*(1) + \frac{dx^*(1)}{dt}, x(1) - \tilde{x}(1) \right\rangle. \tag{39}$$

Finally adding the inequalities (38), (39) we conclude that $J[x(t)] \geq J[\tilde{x}(t)]$, $\forall x(t)$, $t \in [0,1]$ i.e. $\tilde{x}(t)$, $t \in [0,1]$ is optimal. □

Next, we show that if in the problem (PC) a function $W_k (k \in I)$ is only dependent on $x$ and $x'$, the adjoint inclusion involves only one conjugate variable, that is, there is no auxiliary adjoint variable $\psi^*(t)$ in the second order adjoint differential inclusion.

**Corollary 5.1** Assume that for the problem (PC) $W_k(x, x', x'') \equiv W_k^1(x, x'')$ i.e. a function $W_k$ does not depend on the second argument $x'$ and that the conditions of Theorem 5.1 are fulfilled. Then, the sufficient conditions of optimality in Theorem 5.1 are:

(a) $\left( \dfrac{d^2 x^*(t)}{dt^2} - f_x'(\tilde{x}(t), t), -x^*(t) \right) \in \displaystyle\sum_{k \in I(\tilde{x}(t), \tilde{x}''(t))} \alpha_k(t) \partial_{(x, v_2)} W_k^1(\tilde{x}(t), \tilde{x}''(t))$ a.e. $t \in [0,1]$,

(b) $-\dfrac{dx^*(1)}{dt} \in \partial q(\tilde{x}(1))$, $x^*(1) = 0$,

(c) $\alpha_k(t) W_k^1(\tilde{x}(t), \tilde{x}''(t)) = 0$, $\alpha_k(t) \geq 0$ a.e. $t \in [0,1]$.

*Proof.* For the limit case of the conditions in Theorem 4.2, we evidently obtain the conditions (a)-(c) of corollary. Then, in analogy to the proof of Theorem 5.1 we can prove optimality of the trajectory $\tilde{x}(t)$. □

Let us examine an important particular case of the problem (PC) labelled by (P2):

$$\text{minimize } J[x(\cdot)] = \int_0^1 f(x(t), t) dt + q(x(1))$$

(P2) $\quad W_k^2(x(t), x'(t)) \leq 0$, $k \in I$, a.e. $t \in [0,1]$,

$$x(0) = v_0.$$

For the limit case of the conditions in Theorem 4.3, we get the conditions (d)-(f):

(d) $\left( \dfrac{d u^*(t)}{dt} - \mu f_x'(\tilde{x}(t), t), u^*(t) \right) \in \displaystyle\sum_{k \in I(\tilde{x}(t), \tilde{x}''(t))} \alpha_k(t) \partial_{(x, v_1)} W_k^2(\tilde{x}(t), \tilde{x}'(t))$,

(e) $-u^*(1) \in \partial q(\tilde{x}(1))$,

(f) $\alpha_k(t) W_k^2(\tilde{x}(t), \tilde{x}'(t)) = 0$, $\alpha_k(t) \geq 0$.

Show that the conditions (d)-(f) suffice for optimality of $\tilde{x}(t)$.

**Theorem 5.2** Suppose that in the problem (PC) $W_k(x, x', x'') \equiv W_k^2(x, x')$ i.e. a function $W_k, k \in I$ does not depends on the third argument $x''$ and that the conditions in Theorem 5.1 are satisfied. Then, for optimality of the trajectory $\tilde{x}(t)$ in the problem (P2) it is sufficient that there exist functions $u^*(t)$ and $\alpha_k(t) \geq 0$, $k \in I(\tilde{x}(t), \tilde{x}'(t))$ satisfying the conditions (d)-(f).

*Proof* By condition (d) we can write

$$\sum_{k \in I(\tilde{x}(t), \tilde{x}''(t))} \alpha_k(t) W_k^2(x(t), x'(t)) - \sum_{k \in I(\tilde{x}(t), \tilde{x}''(t))} \alpha_k(t) W_k^2(\tilde{x}(t), \tilde{x}'(t))$$

$$\geq \left\langle \frac{du^*(t)}{dt} - \mu f_x'(\tilde{x}(t), t), x(t) - \tilde{x}(t) \right\rangle + \left\langle u^*(t), x'(t) - \tilde{x}'(t) \right\rangle.$$

Then in analogy to the proof the inequality (34),

$$f(x(t), t) - f(\tilde{x}(t), t) \geq \frac{d}{dt} \left\langle u^*(t), x(t) - \tilde{x}(t) \right\rangle.$$

By integrating this inequality over the interval $[0,1]$ and considering that $x(\cdot), \tilde{x}(\cdot)$ are feasible ($x(0) = \tilde{x}(0) = \upsilon_0$), we can write

$$\int_0^1 [f(x(t), t) - f(\tilde{x}(t), t)] dt \geq \left\langle u^*(1), x(1) - \tilde{x}(1) \right\rangle. \tag{40}$$

By the condition (e) of theorem

$$q(x(1)) - q(\tilde{x}(1)) \geq -\left\langle u^*(1), x(1) - \tilde{x}(1) \right\rangle. \tag{41}$$

Finally adding the inequalities (40), (41) we obtain that $J[x(t)] \geq J[\tilde{x}(t)]$, $\forall x(t), t \in [0,1]$ i.e. $\tilde{x}(t), t \in [0,1]$ is optimal. $\square$

**Corollary 5.2** If $W_k, k \in I$ is a continuously differentiable function, the condition (i) in Theorem 5.1 can be simply written as:

$$\frac{d^2}{dt^2} \left[ \sum_{k \in I(\tilde{x}(t), \tilde{x}'(t)\tilde{x}''(t))} \alpha_k(t) (W_k)_{v_2}'(\tilde{x}(t), \tilde{x}'(t), \tilde{x}''(t)) \right]$$

$$- \frac{d}{dt} \left[ \sum_{k \in I(\tilde{x}(t), \tilde{x}'(t)\tilde{x}''(t))} \alpha_k(t) (W_k)_{v_1}'(\tilde{x}(t), \tilde{x}'(t), \tilde{x}''(t)) \right]$$

$$+ \sum_{k \in I(\tilde{x}(t), \tilde{x}'(t)\tilde{x}''(t))} \alpha_k(t) (W_k)_x'(\tilde{x}(t), \tilde{x}'(t), \tilde{x}''(t)) + f_x'(\tilde{x}(t), t) = 0.$$

*Proof.* Indeed, since

$$\partial_{(x, v_1, v_2)} W_k(\tilde{x}(t), \tilde{x}'(t), \tilde{x}''(t)) = \{(W_k)_x'(\tilde{x}(t), \tilde{x}'(t), \tilde{x}''(t)),$$

$$(W_k)_{v_1}'(\tilde{x}(t), \tilde{x}'(t), \tilde{x}''(t)), (W_k)_{v_2}'(\tilde{x}(t), \tilde{x}'(t), \tilde{x}''(t))\}$$

it follows from the condition (i) in Theorem 5.1 that

$$\frac{d^2 x^*(t)}{dt^2} + \frac{d\psi^*(t)}{dt} - f_x'(\tilde{x}(t), t) = \sum_{k \in I(\tilde{x}(t), \tilde{x}'(t)\tilde{x}''(t))} \alpha_k(t) (W_k)_x'(\tilde{x}(t), \tilde{x}'(t), \tilde{x}''(t)),$$

$$\psi^*(t) = \sum_{k \in I(\tilde{x}(t), \tilde{x}''(t))} \alpha_k(t) (W_k)_{v_1}'(\tilde{x}(t), \tilde{x}'(t), \tilde{x}''(t)),$$

$$-x^*(t) = \sum_{k \in I(\tilde{x}(t), \tilde{x}''(t))} \alpha_k(t) (W_k)_{v_2}'(\tilde{x}(t), \tilde{x}'(t), \tilde{x}''(t))$$

which imply the desired result. □

**Corollary 5.3** Suppose that $W_k(x(t), x'(t), x''(t))$, $k \in I$ are linear functions defined as follows $W_k(x, x', x'') = \langle p_k^0, x \rangle + \langle p_k^1, x' \rangle - \langle p_k^2, x'' \rangle - d_k$, where $p_k^i \in \mathbb{R}^n (i = 0, 1, 2)$ and $d_k$ are fixed vectors and numbers, respectively. Then denoting by $P_0, P_1, Q$ the matrices with rows $p_k^0, p_k^1, p_k^2$ respectively, the conditions of Theorem 5.1 consist of the following

(1) $Q^*\lambda''(t) + P_1^*\lambda'(t) - P_0^*\lambda(t) = 0$,

(2) $-x^{*\prime}(1) - P_1^*\lambda(1) \in \partial q(\tilde{x}(1))$, $x^*(1) = 0$,

(3) $\langle P_0\tilde{x}(t) + P_1\tilde{x}'(t) - Q\tilde{x}''(t) - d, \lambda(t) \rangle = 0$, $\lambda(t) \geq 0$

*Proof.* In fact, taking into account for $\lambda(t)$ and $d$ the vector-function with components $\alpha_k(t)$ and $d_k$ respectively, and the fact that $\partial W_k(x, v_1, v_2) = \{p_k^0, p_k^1, p_k^2\}$ we have the desired result. Indeed, in view of the condition (i) of Theorem 5.1 (we remind that $f(x,t) \equiv 0$) we have

$$\frac{d^2 x^*(t)}{dt^2} + \frac{d\psi^*(t)}{dt} = p_k^0 \alpha_k(t), \quad \psi^*(t) = p_k^1 \alpha_k(t), \quad -x^*(t) = p_k^2 \alpha_k(t), \quad k \in I(\tilde{x}(t), \tilde{x}'(t), \tilde{x}''(t)).$$

Obviously, this result coincides with the result of Theorem 4.1 [23] on sufficient optimality conditions for a Mayer problem with second-order polyhedral differential inclusions. □

**Theorem 5.3** Let us consider the "nonconvex" problem (4)-(6) that is $f(\cdot, t): \mathbb{R}^n \to \mathbb{R}^1$ and $W_k(k \in I), q$ are nonconvex functions. Then for optimality of the arc $\tilde{x}(t)$, $t \in [0,1]$ in the problem (PC) it is sufficient that there exist a pair of absolutely continuous functions $\{x^*(t), \psi^*(t)\}$, $t \in [0,1]$, function $\alpha_k(t) \geq 0$, $k \in I(\tilde{x}(t), \tilde{x}'(t), \tilde{x}''(t))$ satisfying the conditions:

($a_1$) $\left( \dfrac{d^2 x^*(t)}{dt^2} + \dfrac{d\psi^*(t)}{dt} - x^*(t), \psi^*(t), -x^*(t) \right)$

$\in \displaystyle\sum_{k \in I(\tilde{x}(t), \tilde{x}'(t), \tilde{x}''(t))} \alpha_k(t) \partial_{(x, v_1, v_2)} W_k(\tilde{x}(t), \tilde{x}'(t), \tilde{x}''(t))$, a.e. $t \in [0,1]$

($b_1$) $f(x,t) - f(\tilde{x}(t), t) \geq \langle x^*(t), x - \tilde{x}(t) \rangle$, $\forall x \in \mathbb{R}^n$,

($c_1$) $q(x) - q(\tilde{x}(1)) \geq -\left\langle \psi^*(1) + \dfrac{dx^*(1)}{dt}, x - \tilde{x}(1) \right\rangle$, $\forall x \in \mathbb{R}^n$, $x^*(1) = 0$,

($d_1$) $\alpha_k(t) W_k(\tilde{x}(t), \tilde{x}'(t), \tilde{x}''(t)) = 0$, $\alpha_k(t) \geq 0$, a.e. $t \in [0,1]$.

*Proof.* By condition ($a_1$) in the nonconvex case

$$\sum_{k \in I(\tilde{x}(t),\tilde{x}'(t),\tilde{x}''(t))} \alpha_k(t) W_k\left(x(t), x'(t), x''(t)\right) - \sum_{k \in I(\tilde{x}(t),\tilde{x}'(t),\tilde{x}''(t))} \alpha_k(t) W_k\left(\tilde{x}(t), \tilde{x}'(t), \tilde{x}''(t)\right)$$

$$\geq \left\langle \frac{d^2 x^*(t)}{dt^2} + \frac{d\psi^*(t)}{dt} - x^*(t), x(t) - \tilde{x}(t) \right\rangle$$

$$+ \left\langle \psi^*(t), \frac{d(x(t) - \tilde{x}(t))}{dt} \right\rangle - \left\langle x^*(t), \frac{d^2(x(t) - \tilde{x}(t))}{dt^2} \right\rangle.$$

This inequality and condition ($b_1$) confirm (34). Hence, the rest of the proof is analogous to the proof of Theorem 5.1. □

**Remark 5.1** We remark that the same optimality problems (PD), (PC) can be investigated for a different boundary value conditions. For example, for simplicity let us consider the problem (PC) with the boundary value $x(0) \in Q_0$, $x(1) \in Q_1$, where $Q_0, Q_1 \subset \mathbb{R}^n$ are convex subsets. For such boundary value problem with a SSDFI it can be shown that the condition (iii) in Theorem 5.1 shall be replaced by the following conditions (presentation of this proof would lead us too far astray from the main themes of this paper and is therefore omitted):

$$-\psi^*(0) - \frac{dx^*(0)}{dt} \in K_{Q_o}^*(\tilde{x}(0)); \quad \psi^*(1) \in K_{Q_1}^*(\tilde{x}(1)), \tag{42}$$

$$-\frac{dx^*(1)}{dt} \in \partial\varphi(\tilde{x}(1)), \quad x^*(0) = x^*(1) = 0.$$

Similar to the proof of Theorem 5.1 using the boundary conditions (42) and the condition $x^*(0) = x^*(1) = 0$, it can be justified (to avoid long calculations we omit it) that for a problem (PC) with SSDFI and the boundary condition $x(0) \in Q_0$, $x(1) \in Q_1$ the relation (38) has the form

$$\int_0^1 [f(x(t),t) - f(\tilde{x}(t),t)] dt \geq \left\langle \psi^*(1) + \frac{dx^*(1)}{dt}, x(1) - \tilde{x}(1) \right\rangle$$

$$- \left\langle \psi^*(0) + \frac{dx^*(0)}{dt}, x(0) - \tilde{x}(0) \right\rangle. \tag{43}$$

Then from the condition $x^{*\prime}(1) \in \partial q(\tilde{x}(1))$ or equivalently from the inequality

$$q(x(1)) - q(\tilde{x}(1)) \geq -\left\langle \frac{dx^*(1)}{dt}, x(1) - \tilde{x}(1) \right\rangle$$

and from the relation (43) we have the desired result. □

**Example 5.1** This numerical example is to show the feasibility and efficiency of the theoretic results obtained:

$$\begin{aligned} & \text{minimize } x(1), \\ & x(t) - 3x'(t) \leq 0, \\ & x(0) = 1, \, t \in [0,1]. \end{aligned} \tag{44}$$

According to problem (P2) $W_1^2(x(t), x'(t)) = x(t) - 3x'(t)$, $f(x(t)) \equiv 0$, $q(x(t)) = x(t)$, where $W_1^2 : \mathbb{R}^2 \to \mathbb{R}^1$ is a linear function, i.e., $W_1^2(x, v_1) = x - 3v_1$. By conditions (*d*)-(*f*) we show that the optimal trajectory is $\tilde{x}(t) = e^{t/3}$ and the minimal value of the problem (44) is $\tilde{x}(1) = e^{1/3}$. First, it can be easily seen that the subdifferential of $W_1^2$ is the gradient vector, $\{1, -3\}$, that is, $\partial_{(x,v_1)} W_1^2(x, v_1) = \{1, -3\}$. By analogy $\partial q(x) = \{1\}$. Then by the condition (d) of sufficient conditions (d)-(e) $\left( \dfrac{d u^*(t)}{dt}, u^*(t) \right) \in \alpha_1(t) \partial_{(x,v_1)} W_1^2(\tilde{x}(t), \tilde{x}'(t))$. Therefore, we have

$$\frac{d u^*(t)}{dt} = \alpha_1(t); \; u^*(t) = -3\alpha_1(t), \; u^*(1) = -1,$$
$$\alpha_1(t)[\tilde{x}(t) - 3\tilde{x}'(t)] = 0, \; \alpha(t) \geq 0. \tag{45}$$

In view of first relations of (45) we have an initial value problem with simple linear differential equation with constant coefficient

$$\frac{d u^*(t)}{dt} + \frac{1}{3} u^*(t) = 0, \quad u^*(1) = -1,$$

the solution of which is $u^*(t) = -e^{\frac{1}{3}(1-t)}$. Then since $u^*(t) = -3\alpha_1(t)$ we have

$$\alpha_1(t) = -\frac{1}{3} u^*(t) = \frac{1}{3} e^{\frac{1}{3}(1-t)}. \tag{46}$$

The nonnegative function $\alpha_1(t) = \dfrac{1}{3} e^{\frac{1}{3}(1-t)}$ in (46) is nonzero over the entire time interval $[0,1]$. Hence, by the second condition of (45) $\tilde{x}(t) - 3\tilde{x}'(t) = 0$ and the solution to the initial value problem

$$\tilde{x}(t) - 3\tilde{x}'(t) = 0, \; x(0) = 1$$

is $\tilde{x}(t) = e^{\frac{1}{3}t}$ and its value at the point $t = 1$ is $e^{1/3}$.

## 6. Conclusion

This work demonstrates a novel discretization approach in order to solve the optimization of Bolza problem with the system of second-order differential inequalities that generally represent diverse mechanisms in engineering. In accordance with the suggested discretization, the problem with the system of second-order discrete-approximation inequalities are examined. As theorems of equivalence for subdifferential inclusions and adjoint discrete and discrete approximate inclusions are fundamental methods when examining optimality conditions for discrete and discrete-approximate problems, this technique has a vital meaning when obtaining second-order adjoint inclusions. Next, necessary and sufficient conditions for optimality in such problems are inferred. Finally, it is seen that the suggested technique can be utilized to solve different optimization problems with higher order discrete and differential inequalities.